\documentclass[twocolumn, switch]{article} 
\usepackage{preprint}
\usepackage{graphicx}
\usepackage{subcaption}
\usepackage{amsmath}  
\usepackage{booktabs} 
\usepackage{color}
\usepackage{subcaption}
\usepackage{amssymb}  

\usepackage{tabularx}

\usepackage[utf8]{inputenc}	
\usepackage[T1]{fontenc}	
\usepackage{xcolor}		
\usepackage[colorlinks = true,
            linkcolor = purple,
            urlcolor  = blue,
            citecolor = cyan,
            anchorcolor = black]{hyperref}	
\usepackage{nicefrac}		
\usepackage{microtype}		
\usepackage{lineno}		
\usepackage{float}			

\usepackage{lipsum}		

\usepackage{newfloat}
\DeclareFloatingEnvironment[name={Supplementary Figure}]{suppfigure}
\usepackage{sidecap}
\sidecaptionvpos{figure}{c}

\usepackage{titlesec}
\titlespacing\section{0pt}{12pt plus 3pt minus 3pt}{1pt plus 1pt minus 1pt}
\titlespacing\subsection{0pt}{10pt plus 3pt minus 3pt}{1pt plus 1pt minus 1pt}
\titlespacing\subsubsection{0pt}{8pt plus 3pt minus 3pt}{1pt plus 1pt minus 1pt}

\usepackage{tikz,xcolor,hyperref}

\definecolor{lime}{HTML}{A6CE39}
\DeclareRobustCommand{\orcidicon}{
	\begin{tikzpicture}
	\draw[lime, fill=lime] (0,0) 
	circle [radius=0.16] 
	node[white] {{\fontfamily{qag}\selectfont \tiny ID}};
	\draw[white, fill=white] (-0.0625,0.095) 
	circle [radius=0.007];
	\end{tikzpicture}
	\hspace{-2mm}
}
\foreach \x in {A, ..., Z}{\expandafter\xdef\csname orcid\x\endcsname{\noexpand\href{https://orcid.org/\csname orcidauthor\x\endcsname}
			{\noexpand\orcidicon}}
}
 
\title{Convergence Filters for Efficient Economic MPC of Non-dissipative Systems}

\usepackage{xwatermark}
\newwatermark[firstpage,color=gray!90,angle=0,scale=0.28, xpos=0in,ypos=-5in]{*correspondence: \texttt{hdfzj@zjut.edu.cn}}

\usepackage{authblk}

\author[1\thanks{\tt{hdfzj@zjut.edu.cn}}]{Defeng He}
\author[1\thanks{\tt{1112103020@zjut.edu.cn}}]{Weiliang Xiong\orcidA{}}
\author[2\thanks{\tt{syli@sjtu.edu.cn}}]{Shaoyuan Li}
\author[3\thanks{\tt{hdu@uow.edu.au}}]{Haiping Du}  

\affil[1]{College of Information Engineering, Zhejiang University of Technology, Hangzhou 310023, PR China}
\affil[2]{Department of Automation, Shanghai Jiao Tong University, and Key Laboratory of System Control and Information Processing, Ministry of Education of China, Shanghai, China} 
\affil[3]{School of Electrical, Computer and Telecommunications Engineering, University of Wollongong, Wollongong NSW 2522, Australia}

\begin{document}
\twocolumn[ 
  \begin{@twocolumnfalse} 
\maketitle
\begin{abstract}
This note presents a novel and efficient Economic Model Predictive Control (EMPC) scheme specifically designed for non-dissipative systems subject to state and input constraints. To address the stability challenge of EMPC for constrained non-dissipative systems, a new concept of convergence filters is introduced. Three alternative convergence filters are designed accordingly to be incorporated into the receding horizon optimization problem of EMPC. To improve online computational efficiency, the variable horizon approach without explicit terminal state constraints is adopted. This design allows for a flexible trade-off among convergence speed, economic performance, and computational burden via simple parameter adjustment. Moreover, sufficient conditions are rigorously derived to guarantee recursive feasibility and stability. The advantages of the proposed EMPC are validated through simulations on a classical non-dissipative continuous stirred-tank reactor. 
\end{abstract}
\vspace{1.5em}
  \end{@twocolumnfalse} 
] 


\section{INTRODUCTION}
Economic Model Predictive Control (EMPC) has attracted significant attention in both academia and industry due to its ability to optimize transient economic performance while ensuring constraint satisfaction. Since the performance index of EMPC is generally not positive definite or convex related to any steady-state points \cite{ellis2014tutorial}, the traditional decreasing arguments in standard stabilizing MPC, e.g., \cite[Th.~2.24]{rawlings2017model}, are unreliable. This often leads to the typical divergence or periodic oscillation behaviors of closed-loop systems under EMPC controllers \cite{rawlings2017model}.

Existing methods for stabilizing EMPCs are dominated by the dissipativity and turnpike principles \cite{angeli10average, ellis2014tutorial, schwenkel2024linearly}. If these conditions are not satisfied, the system can be facilitated to converge by weighting stabilization cost \cite{alamir2017contraction} or state increments \cite{alamir2021new} in the economic cost. However, the related weighted functions and weight values are challenging to determine. An alternative is to impose stability constraints on the Finite Horizon Optimal Control Problem (FHOCP). For instance, in \cite{heidarinejad2012economic}, the stability constraint of a Lyapunov function was utilized to render asymptotic stability. Using the multi-objective idea \cite{he2015stability}, the stabilizing EMPCs were developed in \cite{he2016economic, zavala2015multiobjective} by introducing artificially constructed stabilizing cost functions, which were further extended to event-triggered EMPC \cite{he2023event}. Nevertheless, the heavy online computational burden associated with solving the nonconvex nonlinear FHOCP restricts the application of these methods in fast real-time systems.

Several EMPC algorithms have been developed to enhance the efficiency of online optimization by eliminating explicit terminal state constraints in the FHOCP and dynamically shortening the horizon during the control process. Concerning terminal constraint elimination, Grüne and Stieler \cite{grune2012nmpc, grune2014asymptotic} proposed terminal-free EMPC schemes with practical asymptotic stability based on the turnpike property, which were further extended to non-stationary periodic orbits in \cite{muller2016economic}. Alamir and Pannocchia \cite{alamir2021new} introduced state increment weighting to drive the system toward a neighborhood of the optimal steady state. While terminal-free designs reduce the number of constraints, shortening the horizon can further decrease the dimensionality of decision variables, thereby significantly improving optimization efficiency. For instance, Ellis and Christofides \cite{ellis2014finite} developed a hierarchical EMPC framework based on a pre-specified Lyapunov function, independent of the dissipativity assumption or turnpike property. Nonetheless, constructing a Lyapunov function explicitly can be challenging in practice. In our prior work \cite{xiong2024two}, a terminal-free, two-stage EMPC approach allowing for dynamic horizon adjustment was proposed; however, stability was only guaranteed under non-decreasing horizons. 

This note introduces a novel concept termed the convergence filter, which serves as a constraint to ensure the stability of closed-loop trajectories even if the constrained nonlinear systems are non-dissipative. In contrast to the control barrier function and the safety filter \cite{wabersich2023data}, which restrict the closed-loop trajectories to satisfy the specified constraints, our convergence filter ensures the convergence behavior of the closed-loop trajectories. Building on this innovation, we develop a Variable-Horizon EMPC (VHEMPC) scheme without explicit terminal constraints to significantly enhance online computational efficiency. Through appropriate parameter tuning, the proposed approach enables a flexible trade-off among computational burden, economic performance, and convergence speed, all while preserving rigorous theoretical guarantees of recursive feasibility and stability. The effectiveness and benefits of the method are demonstrated through a case study on a Continuous Stirred Tank Reactor (CSTR) system. 

\noindent
\textit{\textbf{Notation:}} Denote the $n$-dimensional real space as $R^n$. The set of non-negative real numbers is $R_{\geq 0}$. The set of integers not less than $i$ is denoted by $\mathbf{I}_{\geq i}$. We define $\mathbf{I}_{a}^{b}:=\{a,a+1,\dots,b\}$ as the set of integers from $a$ to $b$. The distance of a point $x$ to a compact set $\mathbf{S}$ is defined as $\operatorname{dist}(x,\mathbf{S}):=\min_{y \in \mathbf{S}} \left\| x-y \right\|$. The ceil function is $\left\lceil \cdot \right\rceil$. The $\mathcal{K}$, ${{\mathcal{K}}_{\infty }}$ and $\mathcal{K}\mathcal{L}$-functions refer to the standard definitions \cite{khalil2002nonlinear}.  

\section{PROBLEM DESCRIPTION}
Consider a discrete-time nonlinear system
\begin{equation}
{{x}_{k+1}}=f({{x}_{k}},{{u}_{k}}),\quad k\in {\mathbf{I}_{\ge 0}}
\end{equation}
\noindent
where $k$ is the time index and $f:\text{ }{{R}^{n}}\times {{R}^{m}}\to \text{ }{{R}^{n}}$ is the system dynamics, $x\in {{R}^{n}},u\in {{R}^{m}}$ are system state and control input, respectively. We assume the state $x_k$ is measurable for feedback control. The system (1) is subject to the following constraints: 
\begin{equation}
{{x}_{k}}\in \mathbf{X},\quad {{u}_{k}}\in \mathbf{U},\quad \forall k\in {\mathbf{I}_{\ge 0}}
\end{equation}
where \textbf{X}, \textbf{U} are compact sets. We consider the stage economic cost ${{L}_{e}}:{{R}^{n+m}}\to R$ which may be non-positive with respect to (w.r.t) the optimal steady state $({{x}_{s}},{{u}_{s}})$, determined by solving the following optimization problem: 
\begin{equation}
({{x}_{s}},{{u}_{s}})=\operatorname*{arg\,min}_{(x,u)\in \mathbf{X}\times \mathbf{U}}\ {{L}_{e}}(x,u), \quad \text{s}\text{.t}\text{.}\ x=f(x,u)\text{.}
\end{equation}
We assume that the solution of (3) is well-defined, unique and located at the origin. Otherwise, a suitable coordinate transformation or a tie-breaking rule can be applied. 

The control strategy of this note aims to steer the state to the origin, while efficiently optimizing the economic performance over the time-varying horizon $N_k$, i.e., 
\begin{equation}
{{J}_{e}}({{x}_{k}},{\textbf{u}_{k}},{{N}_{k}})=\sum\nolimits_{i=0}^{{{N}_{k}}-1}{{{L}_{e}}({{x}_{i|k}},{{u}_{i|k}})}
\end{equation}
\noindent
subject to (2) and ${\textbf{u}_k}=\{{u_{0|k}},{u_{1|k}},\dots,{u_{{N_k-1}|k}}\}$. To this end, we deploy a dual-mode MPC strategy that applies a local state feedback control law within a small neighborhood of the origin (to be designed), and an MPC control law outside.

As $L_e$ may not be positive definite, the direct minimization of (4) usually leads to oscillatory behaviors \cite{ellis2014tutorial}. To this end, the following convergence filter is introduced.

\noindent
\textbf{Definition 1.}  A convergence filter of a system ${{x}_{k+1}}=f({{x}_{k}})$ w.r.t the compact set $\mathbf{S} \subseteq {R^n} $ consists of a positive function $\chi :{{R}^{n}}\to R$ and a real sequence ${{\Pi }_{k}}$ such that  
\begin{subequations}
\begin{equation}
\chi (\operatorname{dist}({{x}_{k}},\mathbf{S}))\le {{\Pi }_{k}},\quad \forall k\in {\mathbf{I}_{\ge 0}}
\end{equation}
\begin{equation}
\lim_{k\to \infty } \Pi_{k}=0.
\end{equation}
\end{subequations}
Constraint (5a) prevents the system trajectory from diverging, while (5b) further ensures convergence to the set $\mathbf{S}$. Therefore, we call the pair $(\chi, \Pi_k)$ a convergence filter since it inherently excludes nonconvergent behaviors. Unlike the Lyapunov function, the convergence filter does not require imposing class-$\mathcal{K}$ upper bounds on $\chi$, which may be difficult to construct for general economic costs. Consequently, the filter cannot independently guarantee closed-loop stability; therefore, we design a controller based on the dual-mode paradigm.

\noindent
\textbf{Lemma 1.}  If the system ${x_{k+1}}=f({x_k})$ allows for a convergence filter, then $\lim_{k\to \infty } \operatorname{dist}(x_{k},\mathbf{S})=0$.

Lemma 1 follows directly from the positive definiteness of $\chi$ and the convergence of ${\Pi _k}$. It should be noted, however, that Lemma 1 only concerns the system behavior at infinite time, and constructing a suitable filter for general nonlinear systems is a nontrivial task.

In the next section, we provide a detailed discussion on the VHEMPC algorithm, focusing specifically on how to construct a convergence filter for non-dissipative systems to drive the closed-loop system to the origin.

\section{VHEMPC WITHOUT TERMINAL STATE CONSTRAINTS}
This section first presents the FHOCP formulation and related lemmas; subsequently, we construct three alternative convergence filters.

\subsection{FHOCP and Closed-loop Systems}
Given the horizon $N_k$ at time $k$, the control sequence is computed by the FHOCP, provided as follows: 
\begin{subequations}
\begin{equation}
V_{e}^{*}({{x}_{k}},{{N}_{k}})=\underset{{{\mathbf{u}}_{k}}}{\mathop{\min }}\,{{J}_{e}}({{x}_{k}},{\mathbf{u}_{k}},{{N}_{k}})
\end{equation}
\begin{equation}
\ \text{s}\text{.t}\text{.}\ \ {{x}_{i+1|k}}=f({{x}_{i|k}},{{u}_{i|k}})
\end{equation}
\begin{equation}
{x_{0|k}}={x_k},\quad {u_{i|k}}\in \mathbf{U},\quad {{x}_{i|k}}\in \mathbf{X}, \quad i\in \mathbf{I}_0^{{{N}_{k}}-1}
\end{equation}
\begin{equation}
{{J}_{a}}({{x}_{k}},{\mathbf{u}_{k}},{{N}_{k}})\le \Pi ({{x}_{k}},{{N}_{k}})
\end{equation}
\end{subequations}
in which the optimal control sequence is represented as $\textbf{u}_k^*=\{ u_{0|k}^*,u_{1|k}^*, \dots,u_{{N_k}-1|k}^*\}$ and the corresponding state trajectory is $\textbf{x}_{k}^*=\{ x_{0|k}^*, x_{1|k}^*, \dots, x_{{N_k}|k}^* \}$. Constraint (6d) is a filter-related constraint where the auxiliary cost function $J_a:{R^n}\times {R^{m{N_k}}}\times {\mathbf{{I}}_{\ge 1}}\to {R_{\ge 0}}$ is defined as 
\begin{equation}
{{J}_{a}}\left( {{x}_{k}},{\textbf{u}_{k}},{{N}_{k}} \right)=\sum\nolimits_{i=0}^{{{N}_{k}}-1}{{{L}_{a}}({{x}_{i|k}},{{u}_{i|k}})}+\lambda V({{x}_{{{N}_{k}}|k}})
\end{equation}
where the weighting factor $\lambda \ge 1$ is utilized to enlarge the domain of attraction \cite{limon2006stability}, and ${{L}_{a}}:{{R}^{m+n}}\to R$ is defined to satisfy the following properties:
\begin{equation}
\begin{cases}
  L_{a}(x,u) = \gamma_{0}(x,u), & \forall x \in \mathbf{X}_{f} \\
  L_{a}(x,u) \ge d, & \forall x \notin \mathbf{X}_{f}.
\end{cases}
\end{equation}
where $d >0 $ and ${\gamma_0}:R^{m+n} \to R$ is a positive definite function satisfying ${\gamma_0}(x,u) \ge \underline{\gamma}_0(x)$ for all $x \in \mathbf{X}_f$, $u \in \mathbf{U}$, with $\underline{\gamma}_0(x)$ also being positive definite.  The $V$ and set ${\mathbf{X}_f}$ satisfy the following assumption, which renders local asymptotic stability and serves as a basic axiom in MPC literature \cite{ellis2014tutorial, rawlings2017model}. 

\noindent
\textbf{Assumption 1.}  There exist a local control law $K(\cdot ):{{R}^{n}}\to {{R}^{m}}$, a sublevel set ${\mathbf{X}_{f}}=\{x|V(x)\le \alpha \}\subseteq \mathbf{X},\alpha >0$ and $\mathcal{K}_{\infty}$ functions ${{\alpha }_{1}},{{\alpha }_{2}}$, such that $\forall x\in {\mathbf{X}_{f}}$, we have $K(x) \in \mathbf{U}$ and 
\begin{subequations}
\begin{equation}
{{\alpha }_{1}}(||x||)\le V(x)\le {{\alpha }_{2}}(||x||)
\end{equation}
\begin{equation}
V\left( f(x,K(x)) \right)-V(x)\le - {{\gamma }_{0}}(x,K(x)).
\end{equation}
\end{subequations}

Solving nonlinear programs is typically time-consuming; therefore, we do not impose the additional terminal state constraint $x_{{{N}_{k}}|k}^{*}\in {\mathbf{X}_{f}}$ in (6), which is typically utilized to ensure feasibility and stability, compared to classical multiobjective EMPC \cite{ he2015stability, zavala2015multiobjective, xiong2025learning}. Unlike the stabilizing MPC presented in \cite{limon2006stability}, we cannot guarantee that $x_{{{N}_{k}}|k}^{*}\in {\mathbf{X}_{f}}$ is satisfied implicitly under economic optimization (6). The following Lemma 2 explains the relationship between the optimal prediction trajectory and the terminal set ${\mathbf{X}_{f}}$.

\noindent
\textbf{Lemma 2 \cite{xiong2024two}.} Consider system (1) under Assumption 1, if the FHOCP (6) is feasible with constraint (6d) satisfying
\begin{equation}
\Pi(x_k, N_k) \le N_k d + \lambda \alpha , 
\end{equation}
then, for the optimal predictive trajectory $\mathbf{x}^*_k$, there is
\begin{equation}
\exists i \in \mathbf{I}^N_0 : x^*_{i|k} \in \mathbf{X}_f.
\end{equation}

\textbf{Proof Sketch:} Lemma 2 can be proven by contradiction. $\square$

The equation (11) includes two mutually exclusive cases:

\textbf{Case A.1\text{:}} $x_k = x^*_{0|k} \in \mathbf{X}_f,\; x^*_{i|k} \notin \mathbf{X}_f, \; \forall i \in \mathbf{I}^{N_k}_1$.

\textbf{Case A.2\text{:}} \( \exists i \in \mathbf{I}^{N_k}_1:\; x^*_{i|k} \in \mathbf{X}_f\).

\noindent
for each of which the closed-loop dynamics and horizon $N_{k+1}$ are discussed below.

\textbf{Case A.1\text{:}} The closed-loop system is controlled by $K$ as
\begin{equation}
x_{k+1} = f(x_k, K(x_k)).
\end{equation}

\textbf{Case A.2\text{:}} The system (1) is feedback controlled by the first element of $\textbf{u}^*_k$ as
\begin{equation}
{x_{k+1}}=f({x_k},u_{0|k}^*).
\end{equation}
In both possibilities, the horizon $N_{k+1}$ for the next time step is calculated as 
\begin{equation}
{{N}_{k+1}}=\max \left\{ \left\lceil {{{\tilde{N}}}_{k+1}}{\upsilon_k}+(1-{\upsilon_k}){N_k} \right\rceil +{\varsigma_k},1 \right\}
\end{equation}
where the user-designed sequences ${\upsilon_k}\in [0,1],{\varsigma_k}\in {\mathbf{I}_{\ge 0}},\forall k\in {\mathbf{{I}}_{\ge 0}}$ control horizon increase and decrease, respectively. We emphasize that the calculation of the minimal feasible horizon ${{\tilde{N}}_{k+1}}$ is nontrivial as it is coupled with the specific design of $\Pi ({{x}_{k}},{{N}_{k}})$ in (6d). The details will be discussed in the next subsection. 

Now, to ensure the feasibility at time $k+1$, we define the auxiliary value function $V_{a}^{e}({{x}_{k}},{{N}_{k}})$ as
\begin{equation}
V_a^e({x_k},{N_k})={J_a}({x_k}, \mathbf{u}_{k}^*,{{N}_{k}}).
\end{equation}
The following lemma establishes the existence of a feasible horizon $N_{k+1}$ and an associated feasible solution, which together ensure that the desired decrease condition is satisfied. 

\noindent
\textbf{Lemma 3.}  If Assumption 1 and (10) hold, and FHOCP (6) is feasible at time $k$, then there exists ${{N}_{k\text{+1}}}$ and an associated feasible control sequence ${{\mathbf{\bar{u}}}_{k+1}}({{N}_{k+1}})$ such that 
\begin{subequations}
\begin{equation}
{{J}_{a}}({{x}_{k+1}},{{\mathbf{\bar{u}}}_{k+1}}({{N}_{k+1}}),{{N}_{k+1}})\le V_{a}^{e}({{x}_{k}},{{N}_{k}})-{{L}_{a}}({{x}_{k}},u_{0|k}^{*})
\end{equation}
\begin{equation}
\begin{aligned}
  &{{J}_{a}}({{x}_{k+1}},{{\mathbf{\bar{u}}}_{k+1}}(N_{k+1}^+),N_{k+1}^{+})\\ 
  &\le {{J}_{a}}({{x}_{k+1}},{{\mathbf{\bar{u}}}_{k+1}}({{N}_{k+1}}),{{N}_{k+1}}), \quad
  \forall N_{k+1}^{+} \in {\mathbf{I}_{\ge {N_{k+1}}}}. 
\end{aligned}
\end{equation}
\end{subequations}
\textbf{Proof.} We first separately discuss Case A.1 and Case A.2 to prove (16a) in Parts 1 and 2, and explain (16b) in Part 3. 

\noindent
\textbf{Part 1.} In Case A.1, the condition $x^*_{N_k|k} \notin \mathbf{X}_f$ implies that $V_a^e({x_k},{N_k})\mathop  > \limits^{(15),(7)} \lambda V(x_{{N_k}|k}^*) \geq \lambda \alpha$. The candidate feasible control sequence $\bar{\mathbf{u}}_{k+1}(N_{k+1}), \; \forall N_{k+1} \in \mathbf{I}_{\ge 1}$ is constructed as  
\begin{equation}
\bar{\mathbf{u}}_{k+1}(N_{k+1}) = \{ K(x_{k+1}), K(x_{1|k+1}), \ldots, K(x_{N_{k+1}|k+1}) \}.
\end{equation}
As \( x_k \in \mathbf{X}_f\) under case A.1, and \(\mathbf{X}_f\) is positively invariant for system (1) by Assumption 1, the feasibility of $\bar{\mathbf{u}}_{k+1}(N_{k+1})$
is guaranteed. Moreover, we have
\begin{equation}
  \begin{aligned}
    & J_a(x_{k+1}, \bar{\mathbf{u}}(N_{k+1}), N_{k+1}) \overset{\text{(8),(9b)}}{\mathop{\le }} \lambda V(x_k) - L_a(x_k, u_{0|k}^*) \\
    & \le \lambda \alpha - L_a(x_k, u_{0|k}^*) \le V_a^e(x_k, N_k) - L_a(x_k, u_{0|k}^*).
  \end{aligned}
\end{equation}
Thus, condition (16a) is satisfied.

\textbf{Part 2.} Under Case A.2, we define the time index at which the trajectory first enters the terminal set $\mathbf{X}_f$ as
\begin{equation}
F_{k}^{*}=\underset{i}{\mathop{\min }}\,\{i:x_{i|k}^{*}\in {\mathbf{X}_f}\},\quad i\in \mathbf{I}_{0}^{{{N}_{k}}}.
\end{equation}
For any possible horizon ${{N}_{k}}_{+1}\in {\mathbf{I}_{\ge F_{k}^{*}}}$, we further define that 
\begin{equation}
\wp _{k}^{*}({{N}_{k\text{+1}}})=\underset{i}{\mathop{\max }}\,\{i:x_{i|k}^{*}\in {\mathbf{X}_f}\}, \quad i\in \mathbf{I}_{F_{k}^{*}}^{\min \{{{N}_{k+1}},{{N}_{k}}\}}
\end{equation}
abbreviated as $\wp _{k}^{*}$. The feasible ${{\mathbf{\bar{u}}}_{k+1}}({{N}_{k\text{+1}}})$ is constructed as 
\begin{equation}
\begin{aligned}
{{\mathbf{\bar{u}}}_{k+1}}({{N}_{k\text{+1}}})= &\{u_{1|k}^{*}\dots
 u_{\wp _{k}^{\text{*}}-1|k}^{*}, \\ 
& K({{x}_{\wp_{k}^{\text{*}}-1|k+1}}),\dots,K({{x}_{{{N}_{k+1}}-1|k+1}})\}.
\end{aligned}
\end{equation}
Then, we have
\begin{equation}
\begin{aligned}
  & J_a(x_{k+1}, {{\mathbf{\bar{u}}}_{k+1}}({{N}_{k\text{+1}}}), N_{k+1}) \\
  & \overset{\text{(7)}}{\mathop{=}} \sum\nolimits_{i=0}^{\wp_k^{*}-1} L_a(x_{i|k+1}, u_{i+1|k}^{*})\quad \\
  &+ \left[ \sum\nolimits_{i=\wp_k^{*}}^{N_{k+1}-1} L_a(x_{i|k+1}, K(x_{i|k+1})) + \lambda V(x_{N_{k+1}|k+1}) \right] \\
  & \overset{(\text{9b})}{\mathop{\le}} \sum\nolimits_{i=0}^{\wp_k^{*}-1} L_a(x_{i|k+1}, u_{i+1|k}^{*}) + \lambda V(x_{\wp_k^{*}|k+1}),
\end{aligned}
\end{equation}
Now, we separately discuss two sub-cases based on $\wp _{k}^{*}$. 

\textbf{Sub-case A.2.1:} If $\wp _{k}^{*}={{N}_{k}}$, i.e., $x_{{{N}_{k}}|k}^{*}\in {\mathbf{X}_f}$, we have 
\begin{equation}
\begin{aligned}
  & \sum\nolimits_{i=0}^{\wp _{k}^{*}-1}{{{L}_{a}}({{x}_{i|k+1}},u_{i+1|k}^{\text{*}})}+\lambda V({{x}_{\wp _{k}^{*}|k+1}}) \\ 
 & \le V_{a}^{e}({{x}_{k}},{{N}_{k}})-{{L}_{a}}({x_k},u_{0|k}^{*}) \\ 
 & +\left[ \lambda V({{x}_{\wp _{k}^{*}|k+1}})+{{L}_{a}}({{x}^*_{\wp _{k}^{*}|k}},K({{x}^*_{\wp _{k}^{*}|k}}))-\lambda V({{x}^*_{\wp _{k}^{*}|k}}) \right] \\ 
 & \overset{\text{(9b)}}{\mathop{\le }}\,V_{a}^{e}({{x}_{k}},{{N}_{k}})-{{L}_{a}}({{x}_{k}},u_{0|k}^{*}).
\end{aligned}
\end{equation}
\textbf{Sub-case A.2.2:} If $\wp _{k}^{*}<{{N}_{k}}$, then, we have $x_{{{N}_{k}}|k}^{\text{*}}\notin {\mathbf{{X}}_{f}}$ and ${{x}_{\wp _{k}^{\text{*}}|k+1}} \in {\mathbf{X}_f}$. It follows that 
\begin{equation}
\begin{aligned}
  & \sum\nolimits_{i=0}^{\wp _{k}^{\text{*}}-1}{{{L}_{a}}({{x}_{i|k+1}},u_{i+1|k}^{\text{*}})}+\lambda V({{x}_{\wp _{k}^{\text{*}}|k+1}}) \\
   & \le \sum\nolimits_{i=0}^{\wp _{k}^{\text{*}}-1}{{{L}_{a}}({{x}_{i|k+1}},u_{i+1|k}^{\text{*}})+\lambda V({{x}^*_{{{N}_{k}}|k}})} \\
 & \le V_{a}^{e}({{x}_{k}},{{N}_{k}})-{{L}_{a}}({{x}_{k}},u_{0|k}^{*}) .
\end{aligned}
\end{equation}
\textbf{Part 3.} The tail elements of ${{\mathbf{\bar{u}}}_{k+1}}({{N}_{k\text{+1}}})$, in both (17) and (21), are determined by local control law $K$ applied to the states within $\mathbf{X}_f$. Thus, (16b) can be obtained by iterating (9b). $\square$

In the following subsection, we explain how to construct the sequence $\Pi ({{x}_{k}},{{N}_{k}})$ in (6d) and further provide corresponding iterative processes to calculate ${{\tilde{N}}_{k\text{+1}}}$ in (14).

\subsection{Convergence Filter Construction}
In this subsection, three sequences ${{\Pi }^{i}}({{x}_{k+1}},{{N}_{k+1}}), i\in \mathbf{I}_{1}^{3}$ are artificially constructed. These sequences, along with the adjusted parameter $\kappa \in (0,1]$, are designed to explicitly compromise convergence speed, economic performance, and computational burden. The sequence initialization at time $k=0$ is given by:  
\begin{equation}
\Pi ({{x}_{0}},{{N}_{0}})\le {{N}_{0}}d+\lambda \alpha. 
\end{equation}
Once the optimization problem (6) is solved at time $k$, the corresponding value of $V_{a}^{e}({{x}_{k}},{{N}_{k}})$ can be evaluated according to (15). Then, the first construction method relies on the value function of the stabilizing FHOCP: 
\begin{equation}
{{\Pi }^{1}}({{x}_{k+1}},{{N}_{k+1}})=(1-\kappa )V_{a}^{e}({{x}_{k}},{{N}_{k}})+\kappa V_{a}^{*}({{x}_{k+1}},{{N}_{k+1}})
\end{equation}
where the stable value function $V_{a}^{*}({{x}_{k+1}},{{N}_{k+1}})$ is defined as
\begin{subequations}
\begin{equation}
V_{a}^{*}({{x}_{k}},{{N}_{k}})=\underset{{{\mathbf{u}}_{k}}}{\mathop{\min }}\,{{J}_{a}}({{x}_{k}},{\mathbf{u}_{k}},{{N}_{k}}\text{)}
\end{equation}
\begin{equation}
\ \text{s.t. (6b)-(6c)} .
\end{equation} 
\end{subequations}
The following Iterative Process 1 is used to calculate ${{\tilde{N}}_{k+1}}$. 

\begin{center}
\centering
\begin{tabularx}{\linewidth}{X}
\toprule
\textbf{Iterative Process 1} \\
\midrule
Given the optimal state trajectory $\textbf{x}_k^*$ and state $x_{k+1}$ \\
\textbf{For} $i = F^*_k-1, \dots, N_k$ \\
\hspace{1em} Calculate $ V_a^*(x_{k+1}, i)$ by solving (27) and $\Pi^1(x_{k+1}, i)$ by (26) \\
\hspace{1em} \textbf{If} $V_{a}^{*}({{x}_{k\text{+1}}},i)\le {{\Pi }^{1}}({{x}_{k\text{+1}}},i)$ and ${{\Pi }^{1}}({{x}_{k\text{+1}}},i) \le id+\lambda \alpha$ \\
\hspace{2em} \textbf{Break the loop} \\
\textbf{End for} \\
\textbf{Return} ${{\tilde{N}}_{k+1}}=i$ \\
\bottomrule
\end{tabularx}
\end{center}

Notably, Iterative Process 1 entails solving optimization problem (27) for each possible horizon length $i$. The requirement for solving up to $N_k$ additional online optimization problems (in the worst case) is computationally impractical for systems demanding high efficiency. Therefore, the second construction is given: 
\begin{equation}
\begin{aligned}
{{\Pi }^{2}}({{x}_{k+1}},{{N}_{k+1}}) & =(1-\kappa )V_{a}^{e}({{x}_{k}},{{N}_{k}})+ \\
& \kappa {{J}_{a}}({{x}_{k+1}},{\mathbf{{\bar{u}}}_{k+1}}({{N}_{k+1}}),{{N}_{k+1}})
\end{aligned}
\end{equation}
\noindent
in which ${{\mathbf{\bar{u}}}_{k+1}}({{N}_{k+1}})$ is from (17) or (21), depending on which of Cases \textbf{A.1} or \textbf{A.2} occurs at time $k$. The Iterative Process 2 for ${{\Pi }^{\text{2}}}({{x}_{k+1}},{{N}_{k+1}})$ is provided as follows. 
\begin{center}
\centering
\begin{tabularx}{\linewidth}{X}
\toprule
\textbf{Iterative Process 2} \\
\midrule
Given the optimal state trajectory $\textbf{x}_k^*$ and state $x_{k+1}$ \\
\textbf{For} $i = F^*-1, \dots, N_k$ \\
\hspace{1em} \textbf{If case A.1} occurs, construct ${{\mathbf{\bar{u}}}_{k+1}}(i)$ as in (17) \\
\hspace{1em} \textbf{Else case A.2} occurs, construct ${{\mathbf{\bar{u}}}_{k+1}}(i)$ as in (21)\\
\hspace{2em} Calculate ${{J}_{a}}\left( {{x}_{k+1}},{{{\mathbf{\bar{u}}}}_{k+1}}(i),i \right)$ and ${{\Pi }^{2}}({{x}_{k\text{+1}}},i)$ \\
\hspace{1em} \textbf{If} ${{J}_{a}}\left( {{x}_{k+1}},{{{\mathbf{\bar{u}}}}_{k+\text{1}}}(i),i \right)\le {{\Pi }^{\text{2}}}({{x}_{k\text{+1}}},i)\le id+\lambda \alpha$ \\
\hspace{2em} \textbf{Break the loop} \\
\textbf{End for} \\
\textbf{Return} ${{\tilde{N}}_{k+1}}=i$ \\
\bottomrule
\end{tabularx}
\end{center}

Additionally, a simpler filter sequence without candidate solution construction is defined as: 
\begin{equation}
{{\Pi }^{3}}({{x}_{k+1}},{{N}_{k+1}})=V_{a}^{e}({{x}_{k}},{{N}_{k}})-\kappa {{L}_{a}}({{x}_{k}},u_{0|k}^{*})
\end{equation}
whose corresponding Iterative Process 3 is given as follows.
\begin{center}
\centering
\begin{tabularx}{\linewidth}{X}
\toprule
\textbf{Iterative Process 3} \\
\midrule
Given the optimal state trajectory $\textbf{x}_k^*$ and state $x_{k+1}$ \\
\textbf{For} $i = F^*-1, \dots, N_k$ \\
\hspace{1em} Calculate ${{\Pi }^{3}}({{x}_{k\text{+1}}},i)$ by (29)\\
\hspace{1em} \textbf{If} ${{\Pi }^{3}}({{x}_{k+1}},i)\le id+\lambda \alpha$ \\
\hspace{2em} \textbf{Break the loop} \\
\textbf{End for} \\
\textbf{Return} ${{\tilde{N}}_{k+1}}=i$ \\
\bottomrule
\end{tabularx}
\end{center}

Note that ${{\Pi }^{\text{2}}}({{x}_{k+1}},{{N}_{k+1}})$ and ${{\Pi }^{\text{3}}}({{x}_{k+1}},{{N}_{k+1}})$ do not entail any optimization but only simple function evaluations. Hence, the additional computation is negligible compared with \cite{he2016economic} and \cite{zavala2015multiobjective} under the same horizons. 

At time $k$, once any iterative process is complete, the horizon $N_{k+1}$ is determined by (14), and FHOCP (6) is solved at time $k+1$, rolling forward. Following the dual-mode MPC principle, we define the set $\mathbf{X}_{\psi }:=\{x|V(x)\le \psi \}\subseteq {\mathbf{X}_{f}}, 0<\psi \le \alpha$, where the terminal control law $K(\cdot)$ is activated thereafter. Finally, by choosing any of the iterative processes (Process 1, 2 or 3), we provide the convergence filter-based VHEMPC algorithm below. 

\begin{center}
\centering
\begin{tabularx}{\linewidth}{X}
\toprule
\textbf{VHEMPC with filter convergence} \\
\midrule
\textbf{\textit{OFFLINE}} \\
\hangindent=2em
\textbf{Initialization:} Given $(x_s, u_s)$ and sets $\mathbf{X}_f, \mathbf{X}_{\psi}$, set functions $V$, $L_a$ and $\Pi(x_k, N_k)$, determine $\lambda$, $N_0$ such that (25) holds, choose the sequences ${\upsilon_k},{\varsigma }_{k}$.\\
\textbf{\textit{ONLINE}} \\
\hangindent=2em
\textbf{Step 1.} Measure $x_k$, if $x\in {\mathbf{X}_{\psi }}$, go to Step 2. Otherwise, go to Step 3.\\
\hangindent=2em
\textbf{Step 2.} Control the system by $K(x_k)$, set $k=k+1$. Return to Step 1.\\
\hangindent=2em
\textbf{Step 3.} Solve FHOCP (6) and observe the optimal state trajectory. If Case A.1 occurs, go to Step 4. Otherwise, Case A.2 occurs, go to Step 5.\\
\hangindent=2em
\textbf{Step 4.} Control the system by $K(x_k)$ as in (12). Go to Step 6.\\
\hangindent=2em
\textbf{Step 5.} Control the system by $u_{0|k}^{\text{*}}$ as in (13). Go to Step 6.\\
\hangindent=2em 
\textbf{Step 6.} Compute ${{\tilde{N}}_{k\text{+1}}}$ by chosen Iterative Process and set  horizon as in (14), $k=k+1$. Return to Step 1.\\
\bottomrule
\end{tabularx}
\end{center}

\section{THEORETICAL PROPERTIES OF VHEMPC}
This section analyzes the feasibility and stability properties of the resulting closed-loop system under the designed VHEMPC. We first analyze the relationships among the three convergence filters and prove that ${{\tilde{N}}_{k+1}}$ is well-defined.

\noindent
\textbf{Lemma 4.}  Suppose Assumption 1 holds. If FHOCP (6) is feasible at time $k$ with (10), then 
\begin{equation}
\begin{aligned}
{{\Pi }^{1}}({{x}_{k+1}},{{N}_{k+1}})&  \le {{\Pi }^{2}}({{x}_{k+1}},{{N}_{k+1}}) \\
& \le {{\Pi }^{3}}({{x}_{k+1}},{{N}_{k+1}}), \quad \forall {{N}_{k+1}}\in {\mathbf{{I}}_{\ge {{{\tilde{N}}}_{k+1}}}} 
\end{aligned}
\end{equation}
and ${{\tilde{N}}_{k+1}}\le {{N}_{k}}$ for Iterative Processes 1-3.

\noindent
\textbf{Proof.}  From the construction of three filters, equation (30) can be deduced from the optimality of $V_{a}^{*}({{x}_{k\text{+1}}},{{N}_{k\text{+1}}})$ and the feasibility decrease property (16a).

To prove ${{\tilde{N}}_{k+1}}\le {{N}_{k}}$ holds, we only need to ensure that ${{\tilde{N}}_{k+1}}={{N}_{k}}$ is feasible for all Iterative Processes. From Lemma 3, we have ${{J}_{a}}({{x}_{k+1}},{{\bar{\mathbf{u}}}_{k+1}}({{N}_{k}}),{{N}_{k}})\le V_{a}^{e}({{x}_{k}},{{N}_{k}})$ and hence,
\begin{equation}
{{J}_{a}}\left( {{x}_{k+1}},{{{\bar{\mathbf{u}}}}_{k+\text{1}}}({{N}_{k}}),{{N}_{k}} \right)\overset{\text{(16a),(28)}}{\mathop{\le }}\,{{\Pi }^{\text{2}}}({{x}_{k\text{+1}}},{{N}_{k}})\overset{\text{(10)}}{\mathop{\le }}\,{{N}_{k}}d+\lambda \alpha
\end{equation}
which implies that ${{\tilde{N}}_{k+1}}={{N}_{k}}$ is feasible for Iterative Process 2. Moreover, the optimality of $V_{a}^{*}({{x}_{k\text{+1}}},{{N}_{k}})$ means 
\begin{equation}
V_{a}^{*}({{x}_{k\text{+1}}},{{N}_{k}})\le {{J}_{a}}({{x}_{k+1}},{{\mathbf{\bar{u}}}_{k+1}}({{N}_{k}}),{{N}_{k}})\overset{(16a)}{\mathop{\le }}\,V_{a}^{e}({{x}_{k}},{{N}_{k}}).
\end{equation}
It follows that 
\begin{equation}
\begin{aligned}
V_{a}^{*}({{x}_{k+1}},{{N}_{k}})& \overset{(26)}{\mathop{\le }}\,{{\Pi }^{1}}({{x}_{k+1}},{{N}_{k}})\\
& \overset{(30)} {\mathop{\le }}\,{{\Pi }^{2}}({{x}_{k+1}},{{N}_{k}})\overset{(31)}{\mathop{\le }}\,{{N}_{k}}d+\lambda \alpha .
\end{aligned}
\end{equation}
Thus, ${{\tilde{N}}_{k+1}}={{N}_{k}}$ is feasible for Iterative Process 1. Finally, 
\begin{equation}
{{\Pi }^{3}}({{x}_{k+1}},{{N}_{k}})\overset{(29)}{\mathop{\le }}\,V_{a}^{e}({{x}_{k}},{{N}_{k}})\overset{(10)}{\mathop{\le }}\,{{N}_{k}}d+\lambda \alpha
\end{equation}
means the feasibility of ${{\tilde{N}}_{k+1}}={{N}_{k}}$ for ${{\Pi }^3}({{x}_{k+1}},{{N}_{k}})$. $\square$

Lemma 4 guarantees that the iterative processes terminate no later than ${{\tilde{N}}_{k+1}}={{N}_{k}}$. Moreover, for any returned value of ${{\tilde{N}}_{k+1}}<{{N}_{k}}$, the conditions imposed during the Iterative Processes ensure ${{\Pi }^3}({{x}_{k+1}},{{\tilde{N}}_{k+1}})\le {{\tilde{N}}_{k+1}}d+\lambda \alpha $ holds at any time $k$.

\noindent
\textbf{Remark 1.} Note that a smaller filter sequence generally implies a faster convergence speed. Specifically, we observe that $\underset{\kappa \to 0}{\mathop{\lim }}\,{{\Pi }^{i}}({{x}_{k+1}},{{N}_{k+1}})=V_a^e({x_k},{N_k}), \forall i\in \mathbf{I}_1^3$. Thus, the additional effort involved in optimizing for $V_a^*({x_{k+1}},{N_{k+1}})$ and constructing the candidate sequence ${{\mathbf{\bar{u}}}_{k+1}}({{N}_{k+1}})$ provides greater flexibility to compromise between economic performance and convergence speed. In particular, under the choice of $\kappa =1$ and assuming that (27) has a unique solution, the filter sequence ${{\Pi }^{1}}({x_k},{N_k})$ can restore the closed-loop trajectory of a standard stabilizing MPC minimizing ${J_a}(x_k,{\mathbf{u}_k},{N_k})$.

The core feasibility and the stability properties of the closed-loop system are analyzed as follows.

\noindent
\textbf{Theorem 1.} If Assumption 1 and (25) hold, and the FHOCP (6) is feasible at $k=0$, then, the FHOCP is recursively feasible for ${{x}_{k}}\notin {\mathbf{X}_{\psi}}$ with horizons determined by (14) and ${{\Pi }^{i}}({{x}_{k}},{{N}_{k}})\le {{N}_{k}}d+\lambda \alpha \text{,}$$\forall i\in \mathbf{I}_{\text{1}}^{\text{3}},\forall k\in {\mathbf{{I}}_{\ge 0}}$  holds. 

\noindent
\textbf{Proof.}  We first discuss ${{\Pi }^{2}}({{x}_{k}},{{N}_{k}})$ and ${{\Pi }^{3}}({{x}_{k}},{{N}_{k}})$ in Part 1, then turn to ${{\Pi }^{1}}({{x}_{k}},{{N}_{k}})$ in Part 2.

\noindent
\textbf{Part 1.} We claim that candidate solutions corresponding to Cases A.1 and A.2 have been given in Lemma 3. Constraints (6b) and (6c) are trivial, while the (6d) can be verified as follows: 
\begin{equation}
\begin{aligned}
  & {{J}_{a}}({{x}_{k+1}},{{{\mathbf{\bar{u}}}}_{k+1}}({{N}_{k+1}}),{{N}_{k+1}}) \\
  & \le {{J}_{a}}({{x}_{k+1}},{{{\mathbf{\bar{u}}}}_{k+1}}({{{\tilde{N}}}_{k+1}}),{{{\tilde{N}}}_{k+1}}) \\ 
 & \overset{\text{Process 2}}{\mathop{\le }}\,{{\Pi }^{2}}({{x}_{k+1}},{{{\tilde{N}}}_{k+1}})\overset{(\text{30})}{\mathop{\le }}\,{{\Pi }^{3}}({{x}_{k+1}},{{{\tilde{N}}}_{k+1}}) \\ 
 & \overset{\text{Process 3}}{\mathop{\le }}\,{{{\tilde{N}}}_{k+1}}d+\lambda \alpha \overset{(14)}{\mathop{\le }}\,{{N}_{k+1}}d+\lambda \alpha , \quad \forall {{N}_{k+1}}\in {\mathbf{{I}}_{\ge {{{\tilde{N}}}_{k+1}}}} \\ 
\end{aligned}
\end{equation}
where the first inequality follows the same logic as (16b) and the fact that ${{\tilde{N}}_{k+1}}\le {{N}_{k}}$.  

\noindent
\textbf{Part 2.} For ${{\Pi }^{1}}({{x}_{k}},{{N}_{k}})$, we only need to confirm that (6d) holds for $V_a^*({x_{k+1}}, N_{k+1})$. By iterative process 1, there is  
\begin{equation}
V_{a}^{*}({{x}_{k+1}},{{\tilde{N}}_{k+1}})\le {{\Pi }^{1}}({{x}_{k+1}},{{\tilde{N}}_{k+1}})\le {{\tilde{N}}_{k+1}}d+\lambda \alpha .
\end{equation}
Furthermore, by Proposition 2.18 in \cite{rawlings2017model}, we have 
\begin{equation}
V_{a}^{*}({{x}_{k+1}},{{N}_{k+1}})\le V_{a}^{*}({{x}_{k+1}},{{\tilde{N}}_{k+1}}).
\end{equation}
Therefore, it holds that
\begin{equation}
\begin{aligned}
V_{a}^{*}({{x}_{k+1}},{{N}_{k+1}}) &\le {{\Pi }^{1}}({{x}_{k+1}},{{N}_{k+1}})   \\
&\le {{\tilde{N}}_{k+1}}d+\lambda \alpha \le {{N}_{k+1}}d+\lambda \alpha ,
\end{aligned}
\end{equation}
meaning the FHOCP is feasible at $k+1$.  $\square$

\noindent
\textbf{Theorem 2.} Consider system (1) under the same conditions of Theorem 1, if the initial state ${{x}_{0}}\in \mathbf{X} \backslash {\mathbf{X}_{\psi }}$, then the closed-loop system reaches the set ${\mathbf{X}_{\psi }}$ in at most ${{S}_{\psi }}$ steps, given by  
\begin{subequations}
\begin{equation}
{{S}_{\psi }}={({{N}_{0}}d+\lambda \alpha )}/{\kappa \tau }
\end{equation}
\begin{equation}
\tau = \min_{(x,u) \in \mathbf{X} \times \mathbf{U}} {L_a}(x,u), \quad \text{s.t. } V(x) > \psi \end{equation}
\end{subequations}
with $\tau > 0$. Moreover, the closed-loop system is asymptotically stable thereafter. 

\noindent
\textbf{Proof.} For $V_{a}^{e}({{x}_{k+1}},{{N}_{k+1}})$, we have 
\begin{equation}
\begin{aligned}
  V_{a}^{e}({{x}_{k+1}},{{N}_{k+1}})& \overset{(\text{6}d)}{\mathop{\le }}\,{{\Pi }^{i}}({{x}_{k+1}},{{N}_{k+1}})\\
  & \overset{(30),(29)}{\mathop{\le }}\,V_{a}^{e}({{x}_{k}},{{N}_{k}})-\kappa {{L}_{a}}({{x}_{k}},u_{0|k}^{*}) \\ 
 & \overset{(6d)}{\mathop{\le }}\,{{\Pi }^{i}}({{x}_{k}},{{N}_{k}})-\kappa {{L}_{a}}({{x}_{k}},u_{0|k}^{*})\\
 & \overset{(\text{8})}{\mathop{\le }}\,{{\Pi }^{i}}({{x}_{k}},{{N}_{k}})-\kappa \chi ({{x}_{k}}), \quad \forall i\in \mathbf{I}_{1}^{3}
\end{aligned}
\end{equation}
where $\chi :=\min \{{{\underline{\gamma} }_{0}},d\}$ is positive definite. Thus, ${{\Pi }^{i}}({{x}_{k}},{{N}_{k}})$ is a monotonically decreasing sequence satisfying 
\begin{equation}
\underset{k\to \infty }{\mathop{\lim }}\,{{\Pi }^{i}}({{x}_{k}},{{N}_{k}})=0. 
\end{equation}
Consequently, the pairs $\text{(}\chi ,{{\Pi }^{i}}({x_{k}},{{N}_{k}})),\forall i\in \mathbf{I}_{1}^{3}$ are verified as convergence filters of the closed-loop system w.r.t. the origin. By Lemma 1, the closed-loop system converges to the origin until the state enters ${\mathbf{{X}}_{\psi }}$. Furthermore, (8) and (40) imply that 
\begin{equation}
\begin{aligned}
{{\Pi }^{i}}({{x}_{k+1}},{{N}_{k+1}})\le {{\Pi }^{i}}({{x}_{k}},{{N}_{k}})-\kappa \tau, \quad 
\forall i\in \mathbf{I}_{1}^{3},\forall x\in \mathbf{X} \backslash {\mathbf{X}_{\psi }} . 
\end{aligned}
\end{equation}
Therefore, the upper bound of ${{S}_{\psi }}$ can be derived from (42) and (25). Finally, (9) ensures that the closed-loop system asymptotically stabilizes to the origin.  $\Box$

Theorem 2 implies that the upper bound of infinite-time average economic performance can be obtained: 
\begin{equation}
\mathop {\lim }\limits_{k' \to \infty } {\mkern 1mu} \frac{1}{{k'}}\sum\nolimits_{k = 0}^{k'} {{L_e}({x_k},{u_k})}  \leq {L_e}({x_s},{u_s}).
\end{equation}

\section{DISCUSSION ON PARAMETER DESIGN AND EXTENSIONS}
This section discusses details of parameter and set design, which significantly influence the convergence speed, transient and long-term performance, and online computational burden. We also explore potential extensions to robust, switching and learning techniques.

In the offline phase, the functions $L_a$ and $V$ satisfying (8) and Assumption 1 can typically be chosen as quadratic forms \cite{xiong2025learning, xiong2024two}, and the corresponding sets $\mathbf{X}_f$ and terminal control law $K$ can be computed via linear matrix inequalities  \cite{lazar2018computation}. Subsequently, the parameter $d$ in (8) can be estimated as
\begin{equation}
d = \min_{(x,u) \in \mathbf{X} \times \mathbf{U}} L_a(x,u), \quad \text{s.t. } V(x) > \alpha.
\end{equation}
Note that the value of $d$ obtained through this minimization is generally a conservative lower bound. A smaller $d$ value necessitates a longer initial horizon $N_0$ to satisfy the initialization condition (25), which consequently increases the online computational burden. To reduce $N_0$, one can either increase the value of the weighting factor $\lambda$ or iteratively enlarge the terminal set using advanced reinforcement learning approach in \cite{cui2026characterization}. 

Here, we propose a practical approach to reduce $N_0$ by restricting the Lipschitz constant of $L_a$ outside $\mathbf{X}_f$. For example, we can choose 
\begin{equation}
\begin{cases}
{{L}_{a}}(x,u)={{\gamma }_{0}}(x,u),& x\in {\mathbf{X}_{f}} \\ 
{{L}_{a}}(x,u)=d+b{{\left\| x \right\|}^{2}}, & x\notin {\mathbf{X}_{f}} \\ 
\end{cases}
\end{equation}
where $b$, $d$ are the adjustment parameters. Then, consider any trajectory arriving at $\mathbf{X}_f$ with length $\bar{N}$ and auxiliary cost ${{J}_{a}}({{x}_{0}}, {\textbf{u}_{k}}, \bar{N}\text{)}$. The initial horizon $N_0$ satisfying (25) entails  
\begin{equation}
{{N}_{0}}\ge \frac{{{J}_{a}}({{x}_{0}},{\mathbf{{u}}_{k}},\bar{N}\text{)}-\lambda \alpha }{d}. 
\end{equation}
When $b\to 0,d\to \infty $, we have ${{N}_{0}}=\bar{N}$. Therefore, the length of any trajectory that reaches the terminal set (e.g., a trajectory obtained from a terminal state constraint-based scheme) can be chosen as the initial horizon of VHEMPC. Consequently, a shorter horizon leads to a significant improvement in computational efficiency.


In the first mode, the parameter sequences ${\upsilon_k}$ and ${\varsigma_k}$ directly affect the dimension of the decision variables in (6) via adjusting horizon $N_k$, thereby directly influencing the online computational burden. Note that longer horizon can avoid \textit{myopic behaviours} \cite{ellis2014tutorial} and improves long-term performance, thus, ${\upsilon_k}$ and ${\varsigma_k}$ serve as tuning parameters to compromise long-term performance and computation burden. 

Respecting the choice of convergence filter, a greater sequence ${{\Pi}^{i}}({x_{k}}, {N_{k}})$ allows the system to explore a larger feasible region to search for economically optimal solutions, improving both transient and long-term economic performance. However, this comes at the cost of slower convergence speed. Therefore, the relation in (30) provides guidance to balance these factors. Additionally, once the convergence filter is specified, we can further adjust the value of $\kappa$ to determine ${\Pi^i}({x_k}, {N_k})$, providing additional degrees of freedom for fine-tuning. 

In the second mode inside set ${\mathbf{X}_{\psi}}$, the algorithm only requires $ 0<\psi \le \alpha$. A smaller $\psi$ implies that the system must perform more steps of economic optimization in the first mode, whereas a larger $\psi$ allows earlier application of the terminal controller $K$, reducing computational effort. Finally, we note that the estimate of $S_{\psi}$ is jointly affected by user-designed $\kappa$, $\psi$ and $L_a$, whilst its conservatism can be reduced by increasing $\psi$ and $\kappa$, and by limiting the Lipschitz constant of $L_a$ outside $\mathbf{X}_f$ as in (45). Considering the special case $\psi = \alpha$, $\kappa = 1$, $d \to \infty$, and $b \to 0$, we have ${S_\psi } = ({N_0}d + \lambda \alpha )/\kappa \tau  = ({N_0}d + \lambda \alpha )/d = {N_0}$, whose value has already been discussed above in (46).

Next, we consider several extensions of the proposed VHEMPC. 

\textbf{Robust VHEMPC:} The algorithm can handle bounded additive disturbances $||w|| \leq \bar{w},\, w \in R^n, \bar{w} \geq 0$ using established tube-based or min-max approaches. Our prior work \cite{xiong2025varying} provides input-to-state stability results for increasing horizons. We conjecture that for disturbed systems, the convergence filter should be defined relative to the minimal robust positively invariant set rather than the origin. Correspondingly, the decreasing condition (42) is generalized as: 
\begin{equation}
{\Pi^i}(x_{k+1}, N_{k+1}) \leq {\Pi^i}(x_k, N_k) - \kappa \chi ({{x}_{k}}) + \mu(\bar{w})
\end{equation}
where $\mu$ is a class-$\mathcal{K}$ function. The strict derivation of sufficient conditions for the closed-loop input-to-state stability remains a future topic.  

\textbf{Convergence filter switching:} Theorem 1 implies that the constraint (6d) is satisfied under all three filters for all $k$. Thus, the VHEMPC implementation allows utilizing different filters at different times without violating recursive feasibility. This feature further provides additional flexibility to compromise among computational burden, convergence speed, and transient economic performance. As Theorem 2 holds for any filter, the overall stability results remain valid.

\textbf{Time-varying and unknown stochastic cost:} Note that a main advantage of our VHEMPC is that no dissipativity assumption is required. Since the convergence behavior is ensured by the convergence filters, it is independent of specific economic costs. Consequently, provided that the steady state whose related terminal elements exist is specified, the proposed VHEMPC can handle time-varying and unknown economic costs. For more details, see our recent work \cite{xiong2025varying}, where stochastic unknown costs are captured using a mixed-kernel architecture.

\section{ILLUSTRATIVE EXAMPLE}
A non-dissipative CSTR is employed to demonstrate the effectiveness of the proposed control scheme. The FHOCP is solved using the interior-point algorithm, as implemented in MATLAB R2024b. The system dynamics, cost function, and constraints follow the settings established in \cite{xiong2024two}.

For a fair comparison with the prior work \cite{xiong2024two}, we first adopt the convergence filter ${{\Pi}^{2}}(x_k, N_k)$ and a quadratic auxiliary cost $J_a(x,u)$. The initial state is set as $[0.10, 0.081, 0.095]$. Algorithm parameters are fixed as $\alpha = 0.039$, $\psi = 10^{-20}$, $\lambda = 5$, $\kappa = 0.97$, $N_0 = 57$, and $d = 0.0051$. Four different Parameter Settings (PS) for the horizon adjustment sequences ${\upsilon_k}$ and ${\varsigma_k}$ are tested to evaluate the control effect, and all results are summarized in Fig.1.
\begin{figure}
	\centering
	\includegraphics[width=.5\textwidth]{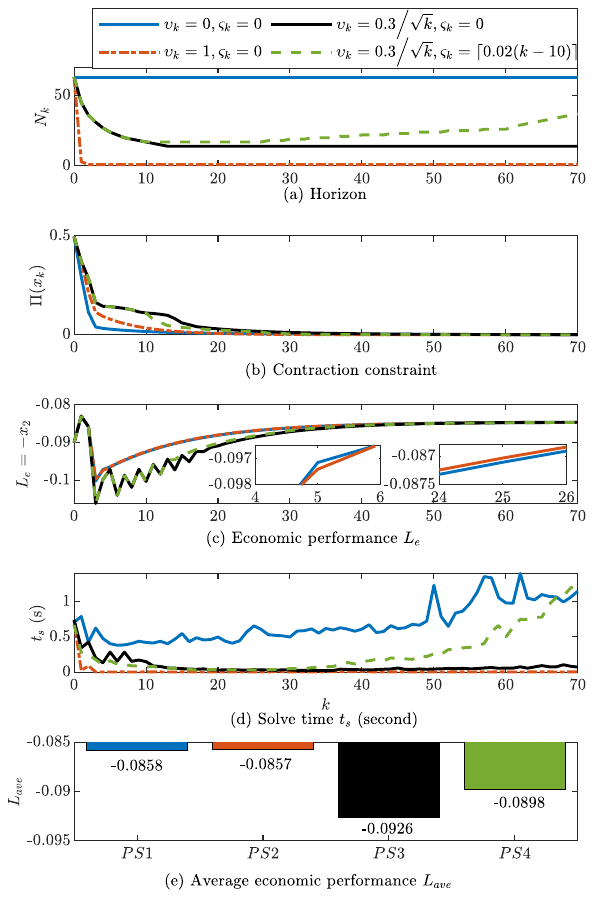}
	\caption{Experimental results of CSTR under different parameters ${{\upsilon }_{k}},{{\varsigma }_{k}}$. (a) Horizon values; (b) Values of convergence filters $\Pi ({{x}_{k}},{{N}_{k}})$; (c) Stage economic cost $L_e$; (d) Online optimization time; (e) Average performance $L_{ave}$ from $k=0$ to $k=70$.}
	\label{FIG:1}
\end{figure}
Fig.1 (a) displays that the horizon can be intuitively adjusted by controlling the parameters, where constant horizons are recovered with ${{\upsilon }_{k}}=0,{{\varsigma }_{k}}=0$ (PS1). In Fig.1 (b), by applying the convergence filter, the control scheme ensures that the sequence $\Pi ({x_k},{N_k})$ stably decreases under admissible horizon change. In contrast, the $\mathrm{VHEMPC}$ scheme proposed in \cite{xiong2024two} may cause $\Pi (x_k, N_k)$ to increase when the horizon decreases, which can lead to severe oscillations in the state trajectory. The corresponding stage cost and optimization time are plotted in Fig.1 (c) and Fig.1 (d), respectively. Due to the conflict between convergence and economic performance, a higher value of $\Pi (x_{k+1},N_{k+1})$ induces better transient performance, while the optimization time and the horizon values have an approximately linear relation. In Fig. 1(c) and (e), we observe that PS3 achieves the best average economic performance $L_{ave}$. This is because its time-varying horizon corresponds to the largest sequence of $\Pi(x_k, N_k)$, as shown in Fig. 1(b), enabling exploration of a wider range of optimal solutions. This is followed closely by PS4. Notably, although PS2 exhibits higher $\Pi(x_k, N_k)$ values and outperforms PS1 at $k = 5$, the longer prediction horizon of PS1 yields slightly better transient behavior throughout the later optimization, resulting in a slightly superior average performance compared with PS2.
 \begin{table}
\centering
\caption{Minimal feasible horizon under different parameters $b$}
\begin{tabular}{ccccccc}
\toprule  
$b$ & 30 & 25 & 20 & 15 & 10 & 5 \\
$N_0$ & 71 & 54 & 37 & 21 & 10 & 4 \\
\bottomrule 
\end{tabular}
\end{table}

To further illustrate that Lipschitz constant limitation can reduce the $N_0$ and improve the computational efficiency, we modify the quadratic $L_a$ to (45) and compute the minimum feasible horizon $N_0$ by (46) for the fixed $x_0$. By setting different parameters $b$, the results are presented in Table 1. The initial horizon decreases significantly as parameter $b$ decreases. Accordingly, the conservative nature of initial horizon estimation is significantly reduced.
\section{CONCLUSION}
This note proposed a novel concept termed \textit{convergence filter}, and incorporated it into EMPC to address the stability issue of the non-dissipative nonlinear systems. We developed a convergence filter-based VHEMPC algorithm with provable feasibility and stability, with three alternative iterative procedures for constructing tailored convergence filters. Notably, the trade-off among online computational burden, convergence speed, and economic performance could be flexibly tuned via user-defined sequence parameters. The merits of the proposed approach were validated through a CSTR simulation. Future work will focus on the extension of this scheme to nonlinear systems subject to bounded or stochastic disturbances and will provide a rigorous quantitative analysis of economic performance under varying prediction horizons.

\small
\bibliographystyle{IEEEtrans}
\bibliography{references} 

\end{document}